\documentclass[11pt,reqno,a4paper]{amsart}

\usepackage[plain]{fullpage}

\usepackage{calrsfs}
\usepackage[OT2,T1]{fontenc}
\usepackage[english]{babel}
\usepackage{newtxtext}
\usepackage{newtxmath}
\usepackage{bbm}
\usepackage{tikz}

\usepackage{etoolbox}
\patchcmd{\section}{\scshape}{\bfseries}{}{}
\makeatletter
\renewcommand{\@secnumfont}{\bfseries}
\makeatother

\usepackage[small, labelfont=bf,up, textfont=up, format=plain]{caption} % control captions in float env

\theoremstyle{plain}
\newtheorem{theorem}{Theorem}%[article]

\newtheorem{lemma}[theorem]{Lemma}
\newtheorem{corollary}[theorem]{Corollary}

\theoremstyle{definition}

\usepackage[unicode=true,pdfusetitle,
 bookmarks=true,bookmarksnumbered=false,bookmarksopen=false,
 breaklinks=true,pdfborder={0 0 0},pdfborderstyle={},backref=false,colorlinks=tr
ue]{hyperref}
\definecolor{myblue}{rgb}{0.09,0.32,0.44} %22-84-113
\hypersetup{pdfborder={0 0 0},pdfborderstyle={},colorlinks=true,linkcolor=myblue,citecolor=myblue,urlcolor=blue}
\usepackage{mathtools}
\usepackage{mathrsfs}

\newcommand{\Z}{\mathbf{Z}}

\newcommand{\Q}{\mathbf{Q}}
\newcommand{\F}{\mathbf{F}}
\renewcommand{\C}{\mathbf{C}}
\newcommand{\mb}{\mathbbm}

\newcommand{\ms}{\mathsf}
\renewcommand{\geq}{\geqslant}
\renewcommand{\leq}{\leqslant}

\DeclareMathOperator{\Res}{Res}

\DeclarePairedDelimiter\abs{\lvert}{\rvert}%
\DeclarePairedDelimiter\norm{\lVert}{\rVert}%

\makeatletter
\let\oldabs\abs
\def\abs{\@ifstar{\oldabs}{\oldabs*}}
\let\oldnorm\norm
\def\norm{\@ifstar{\oldnorm}{\oldnorm*}}
\makeatother

\begin{document}

\date{\today}
\title{Topology of zero sets of polynomials with square discriminant}
\author[D.~Hokken]{David Hokken}
\address{\normalfont Mathematisch Instituut, Universiteit Utrecht, Postbus 80.010, 3508 TA Utrecht, Nederland}
\email{d.p.t.hokken@uu.nl}

\begin{abstract} \noindent
Let $\mathcal{N} \neq \{0\}$ be a fixed set of integers, closed under multiplication, closed under negation, or containing $\{\pm 1\}$. We prove that any zero of a polynomial in $\Z[X]$ whose coefficients lie in $\mathcal{N}$ can be approximated in $\C$ to arbitrary precision by a zero of a polynomial in $\Z[X]$ \emph{with square discriminant} whose coefficients also lie in $\mathcal{N}$. Hence the topology of the closure in $\C$ of the set of zeros of all such polynomials is insensitive to the discriminant being a square, in contrast to the Galois groups of the polynomials.
\end{abstract}

\maketitle

\section{Introduction}

\begin{figure}[t]
\centering
\includegraphics[width=0.98\linewidth]{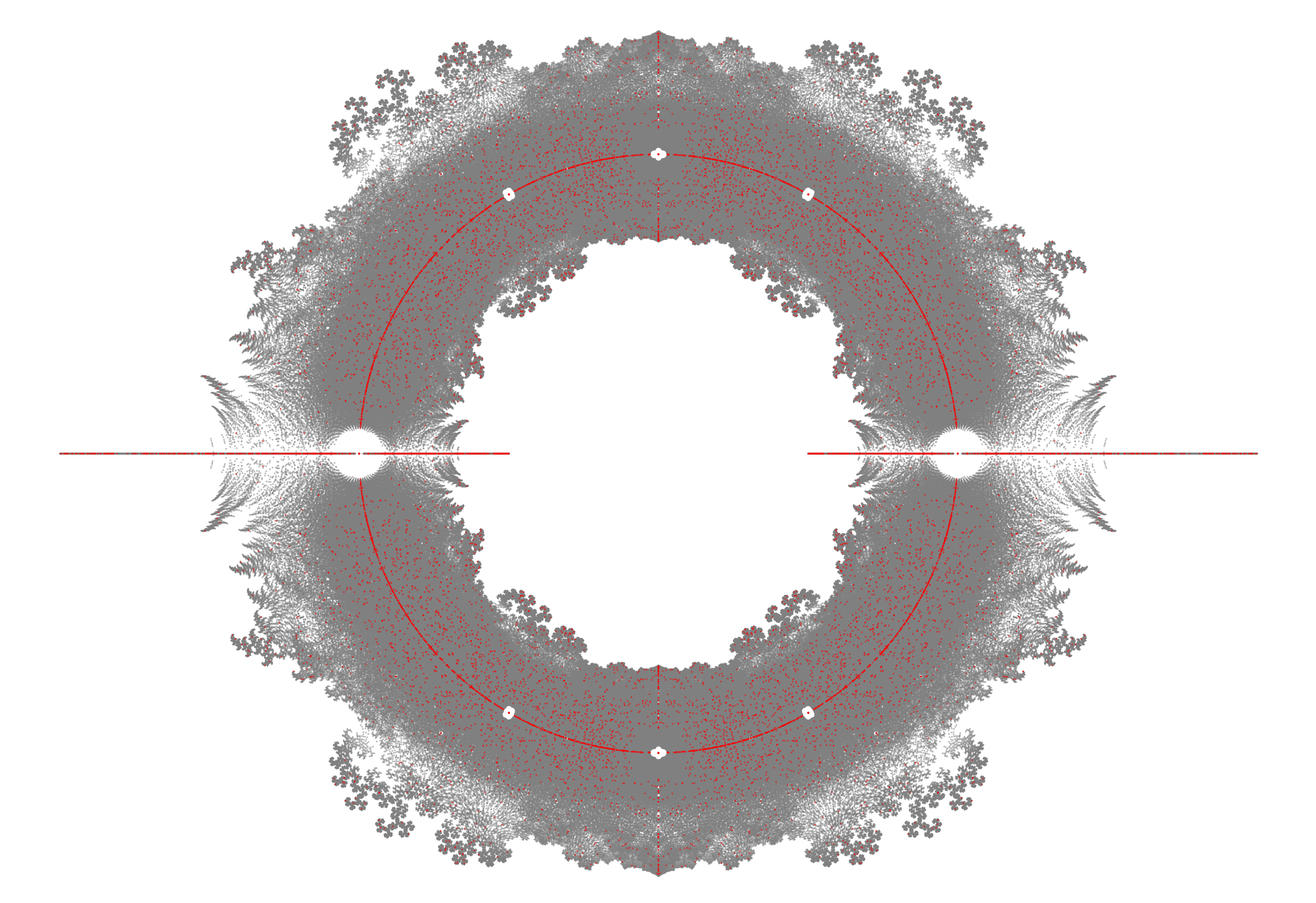}
\caption{Approximation of $M$ (grey) overlain by an approximation of $M^{\square}$ (red) for $\mathcal{N} = \{\pm 1\}$, where only the zeroes of polynomials of degree $n \leq 16$ have been included. The figure is centered at $0 \in \C$, and the `circle' of red dots traces the unit circle. Theorem~\ref{thm:main} asserts that the figure turns red as $n$ tends to infinity.} 
\label{fig:F16}
\end{figure}

Let $\mathcal{N} \neq \{0\}$ be a subset of the integers. Denote by $P = P(\mathcal{N})$ the set of univariate polynomials $f$ whose coefficients lie in $\mathcal{N}$ and such that $f(0) \neq 0$. Define $W=W(\mathcal{N})$ as the countably infinite set of all complex zeroes of all polynomials in $P$ and denote its closure in $\C$ by $M$. First investigated by Barnsley and Harrington \cite{BH} from the point of view of iterated function systems, the set $M$ is a fractal when $\mathcal{N}$ is finite and, as such, abounds in topological features. We mention some of these for the three best-studied cases $M_1 = M(\{0, \pm 1\})$, $M_2 = M(\{\pm 1\})$ (so-called Littlewood polynomials; see Figure~\ref{fig:F16}), and $M_3 = M(\{0, 1\})$. It is a theorem of Bousch \cite{Bousch} that $M_1$ and $M_2$ are connected and locally connected, and of Odlyzko and Poonen \cite{OP} that $M_3$ is path-connected. Bousch also proved the fascinating result that $z \in M_2$ if $z^2 \in M_1$. Bandt \cite{Bandt} showed disconnectedness of the complement $M_1^c$. Subsequent work of Calegari, Koch, and Walker \cite{CKW} established that $M_1^c$ has infinitely many connected components, and also that $M_2^c$ is disconnected. They also confirmed Bandt's conjecture that the interior of $M_1$ is dense away from the real line. 
We refer to \cite{CKW} for further discussion of the related history and results. 

In another, more algebraic direction, one may ask what typical Galois-theoretic properties of the polynomials in $P$ are. For fixed $\mathcal{N}$ of cardinality at least $2$, a folklore conjecture \cite{OP, Konyagin, BV, BKK, BHKP} asserts that a random polynomial $f \in P$ of degree $n$ has the symmetric group $S_n$ as Galois group over $\Q$ with probability tending to $1$ as $n$ tends to infinity. It is an outstanding challenge to rule out the alternating group $A_n$ as likely Galois group. Assume $n>1$, let $a_n$ be the leading coefficient of $f$, and denote the zeroes of $f$ by $\alpha_1, \ldots, \alpha_n \in \C$. Then the Galois group of $f$ is contained in $A_n$ if and only if the discriminant $\Delta(f) \coloneqq a_n^{2n-2} \prod_{i<j} (\alpha_i - \alpha_j)^2$ of $f$ is a nonzero integral square.
Thus the discriminant affects the Galois groups of the polynomials in $P$, and one may ask whether the topology of the respective zero set is also sensitive to the discriminant being a square. 
Let $P^{\square} \subset P$ consist of the polynomials whose discriminant is a square (possibly zero), and define $W^{\square} \subset W$ and $M^{\square} \subset M$ in similar fashion. 
The result here shows that the topology of the closure of the zero set does not change if we restrict to polynomials in $P$ with square discriminant:
\begin{theorem}
\label{thm:main}
Let $\mathcal{N} \neq \{0\}$ be a subset of the integers. Then $M^{\square} = M$ if
\begin{enumerate}
    \item $\mathcal{N}$ is closed under negation, that is $a \in \mathcal{N}$ if and only if $-a \in \mathcal{N}$; or
    \item $\mathcal{N}$ is closed under multiplication, that is $ab \in \mathcal{N}$ if $a, b \in \mathcal{N}$; or 
    \item $\{\pm 1\} \subset \mathcal{N}$.
\end{enumerate}
\end{theorem}
In other words, for any $\mathcal{N}$ as in Theorem~\ref{thm:main}(i), (ii) or (iii), any zero of a polynomial in $P(\mathcal{N})$ can be approximated to arbitrary precision by a zero of a polynomial in $P(\mathcal{N})$ with square discriminant.
Theorem~\ref{thm:main} covers the cases $\mathcal{N} = \{\pm 1\}, \{0,1\}, \{0, \pm 1\}$; one may ask whether $M^{\square} = M$ for any $\mathcal{N} \neq \{0\}$. We cannot exclude the possibility of vanishing discriminant, except in the case $\mathcal{N} = \{\pm 1\}$---see Corollary~\ref{cor:pm1}.

\subsection*{Acknowledgements}
I am grateful to Gunther Cornelissen for discussions on this topic and comments on an earlier version of the paper. 
This research was financially supported by the Dutch Research Council (NWO), project number OCENW.M20.233.

\section{Proofs}
\label{sec:proof}

Let  $f, g \in \Q[X]$ be two nonzero polynomials. Denote by $\Res(f,g) \in \Q$ their resultant. Write $n$ and $a_n$ for the degree and leading coefficient of $f$. If $\alpha_1, \ldots, \alpha_n$ are the complex zeroes of $f$, then the discriminant of $f$ is the rational number $\Delta(f) \coloneqq a_n^{2n-2} \prod_{i<j} (\alpha_i - \alpha_j)^2$. With the convention that the empty product equals $1$, we also adopt this definition when $n \in \{0,1\}$; in particular, the discriminant $\Delta(a_0) = a_0^{-2}$ of a constant $a_0 \in \Q^{\times}$ is a square. Then $\Delta(fg) = \Delta(f)\Delta(g) \Res(f,g)^2$, that is
\begin{equation}
\label{eq:multiplicativity}
\Delta \colon \Q[X] \to \Q/(\Q^{\times})^2 \quad \text{satisfies} \quad \Delta(fg) = \Delta(f)\Delta(g).
\end{equation}
For nonconstant $f$, we also have the formula $\Delta(f) = (-1)^{n(n-1)/2} \Res(f, f')/a_n$.  

We start with the following observation. 

\begin{lemma} 
\label{lem:unitdisk}
The sets $M$ and $M^{\square}$ are equal if and only if they coincide on the open unit disk $\mb{D}$.
\end{lemma}
\begin{proof}
The set $P$ comes equipped with the \emph{reversal map} 
\begin{equation} 
\label{eq:rev}
P \ni f(X) \mapsto f_{\ms{rev}}(X) \coloneqq X^{\deg{f}} f(1/X) \in P
\end{equation}
which is discriminant-preserving by the assumption that $f(0) \neq 0$ if $f \in P$. Hence $\alpha \in W$ if and only if $1/\alpha \in W$. This means that $M$ and $M^{\square}$ are equal if and only if they coincide on the closed unit disk $\overline{\mb{D}}$. For the remainder, it suffices to show that $M^{\square}$ contains the unit circle. 

The zeroes of $p_n(X) \coloneqq 1+X+\cdots+X^n$ lie dense in the unit circle as $n \equiv 1 \bmod{4}$ tends to infinity. We show that $p_n$ has square discriminant, which implies that $ap_n \in P$ has square discriminant for any $a \in \mathcal{N}\setminus\{0\}$, proving the claim.
Applying \eqref{eq:multiplicativity} to $f = p_n$ and $g = X-1$, we find that $\Delta(p_n)$ is a square if
\begin{equation*}
\Delta(X^{n+1}-1) = (-1)^{n(n+1)/2} \Res(X^{n+1}-1, (n+1)X^n) = (-1)^{n(n-1)/2} (n+1)^{n+1}
\end{equation*}
is a square, which holds by the assumption $n \equiv 1 \bmod{4}$. 
\end{proof}

We say $f \in \Q[X]$ is \emph{reciprocal} if it is of even degree and invariant under the reversal map \eqref{eq:rev}.

\begin{lemma}
\label{lem:square disc construction}
Let $f \in \Q[X]$ be of even degree $n$. The following hold:
\begin{enumerate}
    \item If $f$ is reciprocal, then $f$ has square discriminant if and only if $(-1)^{n/2} f(1)f(-1)$ is a square.
    \item The polynomial $f(X)f(X^k)$ has square discriminant for odd $k$.
\end{enumerate}
\end{lemma}
\begin{proof}
For (i), see e.g. \cite[p.~85]{Dubickas}. 
For (ii), if $f \in \Q$, the statement is trivial. We may assume $f$ is monic and does not lie in $\Q$. 
Set $g(X) = f(X^k)$. It suffices to show that $\Delta(f)\Delta(g)$ is a square. Note $g$ has degree $kn$ and $g(\beta) = 0$ if and only if $f(\beta^{k}) = 0$. Hence, with the notation $a_0 = g(0)$,
\begin{equation*}
\Delta(g) = (-1)^{kn(kn-1)/2} \prod_{\beta: g(\beta) = 0} g'(\beta) = (-1)^{kn(kn-1)/2} k^{kn} a_0^{k-1} \prod_{\beta: g(\beta) = 0} f'(\beta^k)
\end{equation*}
where we used $g'(\beta) = f'(\beta^k) k \beta^{k-1}$ in the second equality. On the other hand,
\begin{equation*}
\Delta(f)^k = (-1)^{kn(n-1)/2} \prod_{\alpha: f(\alpha) = 0} f'(\alpha)^k = (-1)^{kn(n-1)/2} \prod_{\beta: g(\beta) = 0} f'(\beta^k).
\end{equation*}
The parity assumptions on $n$ and $k$ imply $\Delta(f)\Delta(g) = (-1)^{n^2k(k-1)/2} k^{kn} a_0^{k-1} \Delta(f)^{k+1}$ is a square.
\end{proof}

Suppose $z \in \mb{D}$. If $F$ is a power series whose coefficients are bounded in absolute value by $H$, then
\begin{equation*}
|F(z)| \leq H(1+|z|+|z|^2+\cdots) = \frac{H}{1-|z|}.
\end{equation*} 
In particular, the radius of convergence of $F$ is at least $1$ and the set of zeroes of $F$ in $\mb{D}$ is well-defined.

\begin{lemma} 
\label{lem:psapprox}
Let $f, g \in \Z\llbracket x\rrbracket$ be two power series whose coefficients are bounded in absolute value. Fix $\epsilon > 0$ and a zero $\alpha \in \mb{D}$ of $f$. Then there is an integer $N$ such that $g$ has a zero $\beta \in \mb{D}$ with $\abs{\alpha-\beta} < \epsilon$ if $X^N \mid (f-g)$. 
\end{lemma}
\begin{proof}
This is an application of Rouch\'e's theorem. Consider an open disk $D$ of radius $R \leq \epsilon$ centered at $\alpha$. Make $R$ small enough such that the closure $\overline{D}$ is contained in the open unit disk and such that $\alpha$ is the only zero of $f$ in $\overline{D}$; this is possible as otherwise the set of zeroes contains an accumulation point in the disk of convergence of $f$, contradicting the identity theorem for analytic functions, see \cite[Theorem~II.3.2]{Lang}. By compactness, $|f|$ attains some positive minimum $m$ on the boundary of $D$. Since the coefficients of $f$ and $g$ are bounded, there is an integer $B$ such that the coefficients of the difference $f-g$ are at most $B$ in absolute value. If the $N$ initial coefficients of $f$ and $g$ coincide, then
\begin{equation*}
|f(z)-g(z)| \leq |z|^N B \sum_{j \geq 0} |z|^j \leq \frac{B |z|^N}{1-|z|}.
\end{equation*} 
Thus choosing $N$ large enough ensures that $|f(z)-g(z)|< m \leq |f(z)|$ for any $z$ lying on the boundary of $D$. By Rouché's theorem \cite[Theorem~VI.1.6]{Lang}, the power series $f$ and $g$ have the same number of zeroes within $D$ (counting multiplicity), so $g$ has a zero $\beta$ satisfying the conditions.
\end{proof}

\begin{proof}[Proof of Theorem \ref{thm:main}]
If $\mathcal{N}$ contains only one element, then the proof of Lemma~\ref{lem:unitdisk} shows that $M = M^{\square}$ is the unit circle. So we may assume the cardinality of $\mathcal{N}$ is at least $2$.
Also by Lemma~\ref{lem:unitdisk}, it suffices to show that for any $\epsilon > 0$ and any complex number $\alpha \in \mb{D}$ that is a zero of some polynomial $f = a_0 + a_1 X + \cdots + a_n X^n \in P$, there is a $\beta \in W^{\square}$ with $\abs{\alpha-\beta} < \epsilon$.
Observe that the power series $F(X) \coloneqq f(X)/(1-X^{n+1})$ also vanishes at $\alpha$. Write $F(X) = b_0 + b_1 X + \cdots$ for its power series expansion. 
For a positive integer $k$ to be selected later, consider the truncation $g(X) = b_0 + \cdots + b_{k-1}X^{k-1}$ of $F$ at the $k$-th term. In case $b_{k-1}=0$, replace $b_{k-1}$ by any nonzero element of $\mathcal{N}$. Then $g \in P$.

\begin{enumerate}
\item If $\mathcal{N}$ is closed under negation, take $k$ even, set $h(X) = g(X) + X^k g(-X)$, and define
\begin{equation*}
f_k(X) = h(X) + aX^{2k} + X^{2k+1} h_{\ms{rev}}(X) \in P
\end{equation*}
where $a \in \mathcal{N}$ can be chosen freely. Then $f_k$ is reciprocal and of degree $4k$ and
\begin{equation*}
f_k(1) = f_k(-1) = g(1)+g(-1)+a-g_{\ms{rev}}(-1)+g_{\ms{rev}}(1),
\end{equation*}
so $f_k(1)f_k(-1)$ is a square. By Lemma~\ref{lem:square disc construction}(i), the polynomial $f_k$ has square discriminant. Furthermore, $F-f_k$ is divisible by $X^{k-1}$. By making $k$ sufficiently large, Lemma~\ref{lem:psapprox} ensures that $f_k$ has a zero $\beta$ with $\abs{\alpha-\beta}<\epsilon$.
\item If $\mathcal{N}$ is closed under multiplication, one can take $k$ odd and $f_k(X) = g(X)g(X^k)$, which lies in $P$ and has square discriminant by Lemma~\ref{lem:square disc construction}(ii). Furthermore, $f_k$ vanishes at all the zeroes of $g$, and $F-g$ is divisible by $X^{k-1}$. By making $k$ sufficiently large, Lemma~\ref{lem:psapprox} ensures that $g$, and hence $f_k$, has a zero $\beta$ with $\abs{\alpha-\beta}<\epsilon$.
\item When $\{\pm 1\} \subset \mathcal{N}$, take $k$ odd, and replace if necessary $b_{k-1}$ by $1$ or $-1$ so as to make $g(1)$ nonzero. Write $\ell=4g(1)-u$ where $u=-1$ if $g(1)>0$ and $u=1$ if $g(1)<0$. Note $u\abs{\ell} = \ell$. Consider the reciprocal polynomial $f_k(X) = (1+X^k)g(X) + uX^{2k}(1+\cdots+X^{\abs{\ell}-1}) + X^{2k+\abs{\ell}}(1+X^k)g_{\ms{rev}}(X)$. Then the degree of $f_k$ is divisible by $4$ and $f_k(1) = f_k(-1) = u$. The same reasoning as in (i) shows that $f_k$ has square discriminant and a zero $\beta$ with $\abs{\alpha-\beta} < \epsilon$ for sufficiently large $k$. 
\end{enumerate}
This finishes the proof.
\end{proof}

From the Galois-theoretic point of view, it would be more interesting to simultaneously control the discriminant and the factorisation of the polynomials $f_k$ constructed in the proof of Theorem~\ref{thm:main}. In (ii), the polynomial $f_k$ is reducible, but in (i) and (iii) one would typically expect $f_k$ to be irreducible. In each construction, one would similarly expect (thanks to the flexibility in choosing $g$) the polynomial $f_k$ to have nonvanishing discriminant, but we cannot ascertain this in general. However, the following result does demonstrate this for the set $\mathcal{N} = \{ \pm 1\}$. To state the result, denote by $M^{\square \neq 0} = M^{\square \neq 0}(\mathcal{N})$ the closure of the set of zeroes of polynomials in $P^{\square}$ with nonvanishing discriminant, and by $C_2 \wr S_{2k}$ the wreath product of the cyclic group $C_2$ with the symmetric group $S_{2k}$.

\begin{corollary}
\label{cor:pm1}
For $\mathcal{N} = \{\pm 1\}$, we have $M = M^{\square \neq 0}$. More precisely, for any $\epsilon > 0$ and zero $\alpha$ of a polynomial with coefficients in $\{\pm 1\}$, there exists a zero $\beta$ of a polynomial with coefficients in $\{\pm 1\}$ and Galois group in $(C_2 \wr S_{4\ell})\cap A_{8\ell}$ for some $\ell > 0$ with the property that $\abs{\alpha-\beta}<\epsilon$.
\end{corollary}
\begin{proof}
In the proof of Theorem~\ref{thm:main}(i), the polynomial $f_k$, with $k>0$ even, is of degree $4k$ and has square discriminant. Over $\F_2$, the polynomial $f_k$ is a factor of the squarefree polynomial $X^{4k+1}-1$, since all its coefficients are $1 \bmod{2}$. Hence $f_k$ is squarefree in $\Q[X]$, which implies that $\Delta(f_k)$ does not vanish, proving that the Galois group $G_k$ of $f_k$ lies in $A_{4k}$; since $f_k$ is reciprocal, in addition $G_k \leq C_2 \wr S_{2k}$. 
\end{proof}

\bibliographystyle{amsplain}

\end{document}